\documentclass{amsart} 

\begin{document}

\newtheorem{theorem}{Theorem}[section]
\newtheorem{lemma}[theorem]{Lemma}
\newtheorem{corollary}[theorem]{Corollary}
\newtheorem{conjecture}[theorem]{Conjecture}
\newtheorem{question}[theorem]{Question}
\newtheorem{problem}[theorem]{Problem}
\newtheorem*{claim}{Claim}
\newtheorem*{criterion}{Criterion}
\newtheorem*{first_thm}{Theorem A}

\theoremstyle{definition}
\newtheorem{definition}[theorem]{Definition}
\newtheorem{construction}[theorem]{Construction}
\newtheorem{notation}[theorem]{Notation}

\theoremstyle{remark}
\newtheorem{remark}[theorem]{Remark}
\newtheorem{example}[theorem]{Example}

\def\area{\text{area}}
\def\id{\text{id}}
\def\cross{\text{cr}} 
\def\H{\mathbb H}
\def\Z{\mathbb Z}
\def\R{\mathbb R}
\def\F{\mathcal F}
\def\out{\textnormal{Out}}
\def\aut{\textnormal{Aut}}
\def\MCG{\textnormal{MCG}}
\def\til{\tilde}
\def\length{\textnormal{length}}
\def\axis{\textnormal{axis}}

\title{Word length in surface groups with characteristic generating sets}
\author{Danny Calegari}
\address{Department of Mathematics \\ Caltech \\
Pasadena CA, 91125}
\email{dannyc@its.caltech.edu}

\date{12/10/2007, Version 0.08}

\begin{abstract}
A subset of a group is characteristic if it is invariant under every
automorphism of the group. We study word length in fundamental groups
of closed hyperbolic surfaces with respect to characteristic generating
sets consisting of a finite union of orbits of the automorphism group,
and show that the translation length of any element with nonzero crossing number
is positive, and bounded below by a constant depending only (and explicitly) on
a bound on the crossing numbers of generating elements. This answers a
question of Benson Farb.
\end{abstract}

\maketitle

\section{Introduction}

Let $S$ be a closed, orientable surface of genus at least $2$.
Then $S$ admits a (nonunique) hyperbolic structure. Let $a \in \pi_1(S)$
and let $[a]$ denote the conjugacy class of $a$. Once we have fixed
a hyperbolic structure on $S$, each nontrivial conjugacy class $[a]$ determines a
unique closed (unparameterized) geodesic $\gamma(a)$ in the corresponding 
free homotopy class on $S$. 

Note that $\gamma(a)$ only depends on the conjugacy class of $a$, so we could
write $\gamma([a])$, but we choose the notation $\gamma(a)$ for simplicity. 

\begin{definition}
If $a$ is primitive, a {\em self-intersection} of $\gamma(a)$ is an unordered
pair of geodesics $\til{\gamma}_1,\til{\gamma}_2$ in the universal
cover of $S$ which cover $\gamma(a)$ and have linked endpoints in the circle
at infinity, up to the action of the deck group $\pi_1(S)$ on such pairs.

The {\em crossing number} of $a$, denoted $\cross(a)$, 
is the number of self-intersections of $\gamma(a)$. If $a = b^n$ then
set $\cross(a) = n^2 \cross(b)$.
\end{definition}

Since $\gamma(a)$ only depends on the conjugacy class of $a$, 
it follows that $\cross$ also depends only on the conjugacy class of $a$.

Moreover, since linking data does not depend on
the choice of hyperbolic structure on $S$ (although the configuration of
$\gamma(a)$ in $S$ typically will), it follows that $\cross(a)$ does not
depend on the choice of hyperbolic structure on $S$.

\begin{lemma}
Crossing number is constant on orbits of $\out(\pi_1(S))$.
\end{lemma}
\begin{proof}
Choose a hyperbolic structure $g$ on $S$, and let $[\phi]$ be an element
of the mapping class group $\MCG(S)$, which
is equal to $\out(\pi_1(S))$ by Dehn-Nielsen, represented by a diffeomorphism 
$\phi\colon S \to S$.
If $[b] = [\phi]([a])$ then $\gamma(b)$ in the $(\phi^{-1})^* g$ metric is the image
of $\gamma(a)$ in the $g$ metric under the diffeomorphism $\phi$.
\end{proof}

For each non-negative integer $n$, let $S_n$ be the subset consisting of all elements 
with crossing number at most $n$, and let $S_n'$ be the subset of $S_n$ consisting
of primitive elements.
 
For each $a \in \pi_1(S)$ define $w_n(a)$ to be the word length of $a$ with respect 
to the generating set $S_n$ and $w_n'(a)$ to be the word length of $a$ with respect 
to the generating set $S_n'$. Obviously $w_n'(a) \ge w_n(a)$ for any $n$, since
$S_n' \subset S_n$.

\begin{remark}
There is an equality $\cross(b) = \cross(b^{-1})$ for every $b \in \pi_1(S)$, so 
the $S_n$ and $S_n'$ are symmetric generating sets.
\end{remark}

\begin{remark}
Each $S_n'$ is a finite union of primitive $\aut(\pi_1(S))$ orbits. Conversely, every
finite union of primitive $\aut(\pi_1(S))$ orbits is contained in $S_n'$ for some $n$, and
every finite union of (not necessarily primitive) $\aut(\pi_1(S))$ orbits is contained in $S_n$ for some $n$.
\end{remark}

Benson Farb asked whether the diameter of $\pi_1(S)$ with respect to $S_0$ is infinite,
and in general to understand the geometry of the Cayley graph of $\pi_1(S)$ with respect
to $\aut(\pi_1(S))$-invariant generating sets. 
The purpose of this note is to prove the following theorem:

\begin{first_thm}
There are constants $C_1(S),C_2(S),C_3(S)$ such that for any non-negative integers $n,m$
and any $a \in \pi_1(S)$ with $\cross(a) > 0$ there is an inequality
$$w_n(a^m) \ge \frac {C_1m} {\sqrt{n} + C_2} - C_3$$
\end{first_thm}

Notice that the constants $C_1,C_2,C_3$ do {\em not} depend on $n$ or on the choice of
$a \in \pi_1(S)$.

\begin{remark}
The dependence of the constants $C_i$ on the surface $S$ is also probably unnecessary.
\end{remark}

\begin{question}
Let $a \in \pi_1(S)$ have $\cross(a)=0$. What is the growth rate of
the word length of $a^m$ with respect to the generating set $S_n'$?
\end{question}

\section{Counting quasimorphisms}

\subsection{Quasimorphisms}

The usual way to obtain lower bounds on word length is to find
a suitable function which is bounded on a generating set, which grows linearly
on powers, and which is almost subadditive under multiplication. A rich source
of such functions is {\em quasimorphisms}.

\begin{definition}
Let $G$ be a group. A quasimorphism is a function $\phi:G \to \R$ for which there
is a smallest real number $D(\phi) \ge 0$ called the {\em defect}, with the
property that for all $a,b \in G$ there is an inequality
$$|\phi(a) + \phi(b) - \phi(ab)| \le D(\phi)$$
A quasimorphism is {\em homogeneous} if $\phi(a^n) = n\phi(a)$ for all integers $n$.
\end{definition}

If $\phi$ is a quasimorphism, the function
$$\overline{\phi}(a):= \lim_{n \to \infty} \frac {\phi(a^n)} n$$
is a homogeneous quasimorphism, whose defect satisfies $D(\overline{\phi}) \le 2D(\phi)$.
Homogeneous quasimorphisms are class functions.
See e.g. \cite{Bavard} or \cite{Calegari_scl} for details.

\subsection{Hyperbolic groups}

We assume the reader is familiar with the basic properties of hyperbolic geometry,
hyperbolic groups and $\delta$-hyperbolic spaces, quasigeodesics, Morse Lemma, 
convexity of distance function, etc. For a reference, see Bridson-Haefliger
\cite{Bridson_Haefliger} or Gromov \cite{Gromov}.

Let $G$ be a group which is $\delta$-hyperbolic with respect to a generating set
$A$. Epstein-Fujiwara \cite{Epstein_Fujiwara}, generalizing a construction due
to Brooks \cite{Brooks}, define so-called {\em counting quasimorphisms} as follows.

\begin{definition}\label{counting_quasimorphism_definition}
Let $\sigma$ be an oriented simplicial path in the Cayley graph $C_A(G)$, and
let $\sigma^{-1}$ denote the same path with the opposite orientation. A
{\em copy} of $\sigma$ is a translate $a\cdot \sigma$ with $a \in G$. For
$\alpha$ an oriented simplicial path in $C_A(G)$, let $|\alpha|_\sigma$ denote
the maximal number of disjoint copies of $\sigma$ contained in $\alpha$. 
For $a \in G$, define
$$c_\sigma(a) = d(\id,a) - \inf_\alpha (\length(\alpha) - |\alpha|_\sigma)$$
where the infimum is taken over all directed paths $\alpha$ in
$C_A(G)$ from $\id$ to $a$.

Define a {\em counting quasimorphism} to be a function of the form
$$h_\sigma(a) := c_\sigma(a) - c_{\sigma^{-1}}(a)$$
\end{definition}

In the sequel we always assume that the length of $\sigma$ is at least $2$. It
is clear from this definition that the homogenization $\overline{h}$ of $h_\sigma$ 
is a class function. A path $\alpha$ as above which achieves the infimum, for
a given $\sigma$ and $a \in G$ is called a {\em realizing path} for $c_\sigma$. 
Since the length of any path is a non-negative integer, a realizing path must exist
for any $a$ and any $\sigma$.

Realizing paths have the following universal property:
\begin{lemma}[Epstein-Fujiwara, Prop. 2.2 \cite{Epstein_Fujiwara}]
Any realizing path for $c_\sigma$ is a $(K,\epsilon)$-quasigeodesic, where
$$K = \frac {\length(\sigma)} {\length(\sigma) - 1}, \quad 
\epsilon = \frac {2\cdot \length(\sigma)} {\length(\sigma) - 1}$$
\end{lemma}

Notice by our hypothesis that $\length(\sigma)$ is at least $2$ that $K \le 2$
and $\epsilon \le 4$. By the Morse Lemma there is a constant $C(\delta)$
such that every realizing path for $c_\sigma$ must be contained in the
$C$-neighborhood of any geodesic from $\id$ to $a$. In particular, one obtains
the following lemma:

\begin{lemma}\label{neighborhood_bound}
There is a constant $C(\delta)$ such that for any path $\sigma \in C_S(G)$
of length at least $2$, and for any $a \in G$, if the $C$-neighborhood of
any geodesic from $\id$ to $a$ does not contain a copy of $\sigma$, then
$c_\sigma(a) = 0$.
\end{lemma}

Finally, the defect of $h_\sigma$ (and therefore of $\overline{h}_\sigma$)
can be controlled independently of $\sigma$:

\begin{lemma}[Epstein-Fujiwara, Prop. 2.13 \cite{Epstein_Fujiwara}]\label{defect_bound}
Let $\sigma$ be a path in $C_S(G)$ of length at least $2$. Then there is
a constant $C(\delta)$ such that $D(h_\sigma) \le C$.
\end{lemma}

\section{Stability of crossings}

If two geodesics in the hyperbolic plane intersect with a small angle, 
then as the geodesics are moved around slightly, 
the point of intersection might move wildly. Nevertheless, two geodesics which cross
are close only in a compact region; for a small perturbation, the point of intersection 
cannot move outside that compact region.

Let $\alpha,\beta$ be two complete geodesics in the hyperbolic plane which cross transversely
at $p$. Fix a small $\epsilon>0$. Let $\alpha_\epsilon$ denote the subset of $\alpha$
which intersects the $2\epsilon$-neighborhood of $\beta$, and similarly define
$\beta_\epsilon$. Let $\alpha',\beta'$ be two geodesics which contain segments
$\alpha'_\epsilon,\beta'_\epsilon$ which are $\epsilon$-close to $\alpha_\epsilon,\beta_\epsilon$
respectively. Then $\alpha',\beta'$ have a transverse intersection which is contained in
$\alpha'_\epsilon \cap \beta'_\epsilon$. We call such a pair of segments $\alpha_\epsilon,\beta_\epsilon$
an {\em $\epsilon$-trap} for the intersection $p$. Or, if $\epsilon$ is understood, just a
{\em trap} for $p$.

Now let $S$ be a closed hyperbolic surface, and let $\gamma \subset S$ be a closed geodesic.
Suppose $\alpha,\beta$ project to $\gamma$ by the covering projection.

\begin{lemma}\label{trap_double_embedded}
Let $\alpha,\beta,\gamma,S$ be as above, and let $\alpha_\epsilon,\beta_\epsilon$ be an
$\epsilon$-trap for $p$, which projects to a self-intersection of $\gamma$. 
Suppose $8\epsilon$ is less than the length of the shortest nontrivial curve on $S$. Then
$\length(\alpha_\epsilon)< \length(\gamma) + 4\epsilon$
and similarly for $\beta_\epsilon$.
\end{lemma}
\begin{proof}
By the definition of an $\epsilon$-trap, the endpoints of $\alpha_\epsilon$ and
$\beta_\epsilon$ are $2\epsilon$ apart. By
the triangle inequality, there is an estimate
$$|\length(\alpha_\epsilon) - \length(\beta_\epsilon)|<4\epsilon$$
Suppose $\length(\alpha_\epsilon) \ge \length(\gamma)+4\epsilon$ and therefore
$\length(\beta_\epsilon) \ge \length(\gamma)$. By hypothesis there are elements
$a,b \in \pi_1(S)$ in the conjugacy class of $\gamma$ which stabilize $\alpha,\beta$
respectively and act as translations through a distance $\length(\gamma)$. Let $q$ be
an endpoint of $\alpha_\epsilon$. Then $d(b^{-1}a(q),q) \le 8\epsilon$. But by
hypothesis, this is shorter than the length of the shortest nontrivial curve on $S$,
so $b=a$ and $\beta = \alpha$, contrary to hypothesis. It follows that
$\length(\alpha_\epsilon) < \length(\gamma) + 4\epsilon$ and the lemma is proved.
\end{proof}

\begin{lemma}\label{stable_intersection}
Let $S,\epsilon$ be as above.
Let $\gamma$ be a closed geodesic on $S$ with a transverse self-intersection. Let
$\alpha$ be a geodesic in the hyperbolic plane covering $\gamma$, and let $\sigma$ be a
geodesic segment in the hyperbolic plane which is $\epsilon$-close to a segment of
$\alpha$, and satisfies $\length(\sigma) > 2\cdot\length(\gamma)+4\epsilon$.
Then the projection of $\sigma$ to $S$ has a transverse self-intersection.
\end{lemma}
\begin{proof}
Let $\beta$ be another geodesic in the hyperbolic plane covering $\gamma$ such
that $\alpha \cap \beta$ contains a point $p$ projecting to a transverse self-intersection
of $\gamma$, and let $\alpha_\epsilon,\beta_\epsilon$ be an $\epsilon$-trap for $p$.
By Lemma~\ref{trap_double_embedded}, $\length(\alpha_\epsilon) + \length(\gamma) \le \length(\sigma)$
so there is a translate of $\sigma$ by an element of the deck group which contains a
segment $\alpha'_\epsilon$ which is $\epsilon$-close to $\alpha_\epsilon$. Similarly,
there is a translate of $\sigma$ containing a segment $\beta'_\epsilon$ which is
$\epsilon$-close to $\beta_\epsilon$. Since $\alpha_\epsilon,\beta_\epsilon$ is an
$\epsilon$-trap for $p$, the segments $\alpha'_\epsilon, \beta'_\epsilon$ intersect
transversely. The projections of $\alpha'_\epsilon$ and $\beta'_\epsilon$ are
both contained in the projection of $\sigma$, and therefore this projection contains
a point of self-intersection, as claimed.
\end{proof}

\begin{remark}\label{Leininger_remark}
Chris Leininger has observed that for every non-simple primitive homotopy class of loop on
a surface $S$, there is some hyperbolic structure on $S$ for which the geodesic representative
intersects itself at a definite angle. This observation, together with the
fact that $w_n$ is characteristic, could be used in place
of the results in this section in the proof of Theorem A.
\end{remark}

\section{Proof of Theorem A}

We now give the proof of Theorem A.

\begin{proof}
Fix a hyperbolic structure on $S$.
Fix a generating set $A$ for $\pi_1(S)$, 
and let $\delta$ be such that
the Cayley graph with respect to this generating set is $\delta$-hyperbolic.
Let $K,\epsilon'$ be such that the Cayley graph $C_A(\pi_1(S))$
is $(K,\epsilon')$ quasi-isometric
to the hyperbolic plane (identified with the universal cover of $S$) 
by a fixed equivariant quasi-isometry. Note that $K,\epsilon'$ can be chosen to
depend only on $S$.

Let $a \in \pi_1(S)$ have $\cross(a)>0$. 
There is a constant $C_1(\delta,|A|)$ depending only on $\delta$ and the cardinality $|A|$
such that $a^{C_1}$ has an axis in $C_A(\pi_1(S))$. We replace $a$ by a suitable
conjugate $b$ of $a^{C_1}$ whose axis contains $\id$. Note that $\gamma(b)$ is
just a $C_1$-multiple of $\gamma(a)$, and therefore contains $C_1^2\cross(a)$ transverse self-intersections.
Since $C_1$ depends only on $\delta,|A|$, and therefore only on $S$, 
it suffices to prove the theorem with $b$ in place of $a$.

\vskip 12pt

Let $l \subset C_A(\pi_1(S))$ be the (oriented) axis of $b$, which by hypothesis
passes through $\id$. Let $N$ be a sufficiently large number (we will say how large
in a moment), and let $\sigma$ be the oriented segment of $l$ from $\id$ to
$b^N$. Define the quasimorphism $h_\sigma$ as in 
Definition~\ref{counting_quasimorphism_definition}.

Claim 1 of Theorem A$^\prime$ in \cite{Calegari_Fujiwara} says that if the translation
length of $b^N$ is sufficiently long (depending only on $\delta$ and $|A|$), 
then $c_{\sigma^{-1}}(b^{Nm}) = 0$ for all positive $m$. From the definition,
$c_\sigma(b^{Nm})=m$ for all positive $m$,
and therefore we have an equality $h_\sigma(b^{Nm}) = m$, valid (by symmetry) 
for all integers $m$.
We can make the translation length of $b^N$ long enough 
by replacing $b$ by a proper power if necessary,
absorbing the constant (multiplicatively) into $C_1$ as above.

On the other hand, let $e \in S_n$. Observe that $e^{C_1}$ has an axis $l_e$.
Note that $\cross(e^{C_1}) \le {C_1}^2 n$ by the definition of $S_n$. 
How many copies of $\sigma$ can there be in a realizing path
for $c_\sigma$ on $e^{C_1}$? Each copy of $\sigma$ in such a realizing path defines
a segment in $l_e$ which is in the $C_2$-neighborhood of a translate of $l$,
where $C_2$ is as in Lemma~\ref{neighborhood_bound}.

It follows that there is a constant $C_3$, depending only on $C_2, K,\epsilon'$,
such that suitable geodesics $\til{\gamma}(e)$ and $\til{\gamma}(b)$ covering $\gamma(e),\gamma(b)$
respectively are distance at most $C_3$ apart on segments of length at least $N\cdot\length(\gamma(b))$.

\vskip 12pt

By convexity of distance in hyperbolic space,
for any $\epsilon > 0$ there is a constant $C_4(\epsilon)$ such that these lifts are
$\epsilon$-close on segments on length at least $N\cdot\length(\gamma(b)) - C_4$.
This tells us how to choose $N$. We choose $\epsilon$ such that $8\epsilon$ is less than
the length of the shortest nontrivial loop on $S$, and choose $N$ so that
$$N\cdot\length(\gamma(b)) - C_4 > 2\cdot\length(\gamma(b)) + 4\epsilon$$
in order to be able to apply Lemma~\ref{stable_intersection}.

Note that $\length(\gamma(b)) \ge 8\epsilon$, and $\epsilon$ 
is bounded below by a positive constant depending
only on the geometry of $S$, so $N,\epsilon$ as above depend only on $S$.

Let $\tau$ be a segment of $\til{\gamma}(e)$ coming from a copy of $\sigma$ in
a realizing path for $c_\sigma$ on $e^{C_1}$, which is
$\epsilon$-close to a segment of $\til{\gamma}(b)$. 
Since $\tau$ by hypothesis has length at least $2\cdot\length(\gamma(b))+4\epsilon$,
by Lemma~\ref{stable_intersection}, the projection of $\tau$ to
$S$ contributes at least one transverse self-intersection to $\gamma(e^{C_1})$.

\vskip 12pt

If a realizing path for $e^{C_1}$ contains $p$ disjoint copies
of $\sigma$, we obtain $p$ such segments $\tau_1,\dots,\tau_p$, and therefore at least
$p^2$ self-intersections of $\gamma(e^{C_1})$.
Since $\cross(e^{C_1}) = C_1^2 \cross(e) \le C_1^2 n$ we obtain an estimate $p^2 \le C_1^2 n$.

Obviously a similar estimate holds for copies of $\sigma^{-1}$ in
$e^{C_1}$. Therefore we obtain
$$|h_\sigma(e^{C_1})| \le p \le C_1 \sqrt{n}$$
Note that after homogenizing, we will obtain a similar estimate for $e$.

\vskip 12pt

By Lemma~\ref{defect_bound}, the defect of $h_\sigma$ depends only on $S$, so
we homogenize $h_\sigma$ and get a homogeneous quasimorphism $\overline{h}_\sigma$ 
with defect depending only on $S$, and satisfying
$$\overline{h}_\sigma(b^m) = \frac m N, \quad 
\overline{h}_\sigma(e) \le C_5\sqrt{n} + C_6$$
for any integer $m$ and any $e \in S_n$, where $N,C_5,C_6$ depend only on $S$.

By the defining property of quasimorphisms (that they are additive with bounded error),
we get a bound of the form
$$w_n(b^m) \ge \frac {C_7m} {\sqrt{n} + C_8} - C_{9}$$
where all constants depend only on $S$.
Since $w_n(a^{mC_1}) = w_n(b^m)$, and $a$ is an arbitrary element with $\cross(a)>0$,
this proves Theorem A.
\end{proof}

\begin{remark}
In view of the estimates in this section, one might be tempted to guess that crossing
number itself is an appropriate function with which to measure word length in
$\pi_1(S)$. Nevertheless, one can construct examples of pairs of elements
$a,b \in \pi_1(S)$ such that $\cross(a) = \cross(b) = 0$ but $\cross(ab)$ is as large as desired.
For example, let $\gamma(b)$ be a very long simple geodesic which is very close
(in the Hausdorff topology) to a minimal geodesic lamination, and let $\gamma(a)$
cross a collection of $10^{100}$ almost parallel strands of $\gamma(b)$.
Let $p \in \gamma(a) \cap \gamma(b)$ 
be a basepoint, in order to pin down the based homotopy classes of $a,b$. Then the
conjugacy class of $ab$ contains a representative loop (not a geodesic) which consists of
the geodesic $\gamma(a)$ followed by the geodesic $\gamma(b)$. Most of the self-intersections
of the loop $ab$ are essential, and therefore the geodesic $\gamma(ab)$ also has at least
$10^{100}$ transverse self-intersections.
\end{remark}

\section{Acknowledgments}

While writing this paper I was partially supported by NSF grant DMS 0405491.
I would like to thank Benson Farb for posing the question which this
paper addresses, Chris Leininger for the observation mentioned in
Remark~\ref{Leininger_remark}, and the anonymous referee. 
This research was carried out while I was at the
Tokyo Institute of Technology, hosted by Sadayoshi Kojima. I thank TIT
and Professor Kojima for their hospitality.

\bibliographystyle{amsplain}

\end{document}